\documentclass{amsart}
\usepackage{geometry,amsmath,amssymb,graphicx,bbm}
\renewcommand{\email}[1]{{\textit{Email:} \texttt{#1}}}
\newcommand{\mathd}{\mathrm{d}}
\newcommand{\nin}{\not\in}
\newcommand{\tmop}[1]{\ensuremath{\operatorname{#1}}}
\newcommand{\tmtextit}[1]{{\itshape{#1}}}
\begin{document}

\title{A Finite Reflection Formula For A Polynomial Approximation To The
Riemann Zeta Function}\author{Stephen
Crowley}\thanks{\email{stephen.crowley@mavs.uta.edu}}\maketitle

\begin{abstract}
  The Riemann zeta function can be written as the Mellin transform of the unit
  interval map $w \left( x \right) = \left\lfloor x^{- 1} \right\rfloor \left(
  x \left\lfloor x^{- 1} \right\rfloor + x - 1 \right)$ multiplied by $s
  \frac{s + 1}{s - 1}$. A finite-sum approximation to $\zeta \left( s \right)$
  denoted by $\zeta_w \left( N ; s \right)$ which has real roots at $s = - 1$
  and $s = 0$ is examined and an associated function $\chi \left( N ; s
  \right)$ is found which solves the reflection formula $\zeta_w \left( N ; 1
  - s \right) = \chi \left( N ; s \right) \zeta_w \left( N ; s \right)$. A
  closed-form expression for the integral of $\zeta_w \left( N ; s \right)$
  over the interval $s = - 1 \ldots 0$ is given. The function $\chi \left( N ;
  s \right)$ is singular at $s = 0$ and the residue at this point changes sign
  from negative to positive between the values of $N = 176$ and $N = 177$.
  Some rather elegant graphs of $\zeta_w \left( N ; s \right)$ and the
  reflection functions $\chi \left( N ; s \right)$ are also provided. The
  values $\zeta_w \left( N ; 1 - n \right)$ for integer values of $n$ are
  found to be related to the Bernoulli numbers. 
\end{abstract}

{\tableofcontents}

\section{The Riemann Zeta Function as the Mellin Transform of a Unit Interval
Map}

The Riemann zeta function can be written as the Mellin transform of the unit
interval map $w \left( x \right) = \left\lfloor x^{- 1} \right\rfloor \left( x
\left\lfloor x^{- 1} \right\rfloor + x - 1 \right)$ multiplied by $s \frac{s +
1}{s - 1}$. {\cite{two-new-zeta-constants}}{\cite{ithsm}}

\begin{figure}[h]
  \resizebox{12cm}{12cm}{\includegraphics{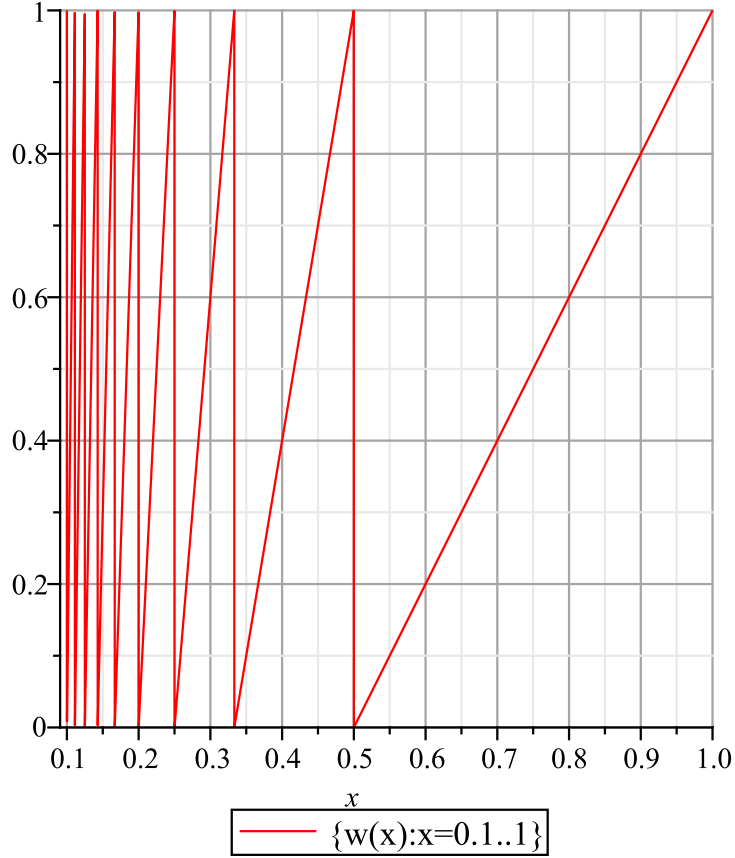}}
  \caption{\label{wfig}The Harmonic Sawtooth map}
\end{figure}
\begin{equation}
  \begin{array}{lll}
    \zeta_w (s) & = \zeta (s) \forall - s \nin \mathbbm{N}^{\ast} & \\
    & = s \frac{s^{} + 1}{s - 1} \int_0^1 \left\lfloor x^{- 1} \right\rfloor
    \left( x \left\lfloor x^{- 1} \right\rfloor + x - 1 \right) x^{s - 1}
    \mathd x & \\
    & = s \frac{s^{} + 1}{s - 1} \sum_{n = 1}^{\infty} \int_{\frac{1}{n +
    1}}^{\frac{1}{n}} n (x n + x - 1) x^{s - 1} \mathd x & \text{ }\\
    & = \sum_{n = 1}^{\infty} s \frac{s^{} + 1}{s - 1} \left( - \frac{n^{1 -
    s} - n (n + 1)^{- s} - s n^{- s}}{s \left( s + 1 \right)} \right) & \\
    & = \sum_{n = 1}^{\infty} \frac{n (n + 1)^{- s} - n^{1 - s} + s n^{-
    s}}{s - 1} & \\
    & = \frac{1}{s - 1} \sum_{n = 1}^{\infty} n (n + 1)^{- s} - n^{1 - s} + s
    n^{- s} & 
  \end{array} \label{sawzeta}
\end{equation}

\subsection{The Truncated Zeta Function}

The substition $\infty \rightarrow N$ is made in the infinite sum appearing
the expression for $\zeta_w (s$) to get a finite polynomial approximation
\begin{equation}
  \begin{array}{ll}
    \zeta_{w_{}} (N ; s) & = \frac{1}{s - 1} \sum_{n = 1}^N n (n + 1)^{- s} -
    n^{1 - s} + s n^{- s}\\
    & = \frac{1}{s - 1} \left( s + (N + 1)^{1 - s} - 1 + s \sum_{n = 2}^N
    n^{- s} - \sum_{n = 2}^{N + 1} n^{- s} \right)\\
    & = \frac{N}{\left( s - 1 \right) \left( N + 1 \right)^s} - \frac{\cos
    \left( \pi s \right) \Psi \left( s - 1, N + 1 \right)}{\Gamma \left( s
    \right)} + \zeta \left( s \right) \forall s \in \mathbbm{N}^{\ast}
  \end{array}
\end{equation}
with equality in the limit except at the negative integers
\begin{equation}
  \begin{array}{ll}
    \lim_{N \rightarrow \infty} \zeta_{w_{}} (N ; s) & = \zeta (s) \forall - s
    \nin \mathbbm{N}^{\ast}
  \end{array}
\end{equation}
and where $\Psi \left( x, n \right) = \frac{\mathd}{\mathd x^n} \Psi \left( x
\right)$ is the polygamma function and $\Psi \left( x \right) =
\frac{\mathd}{\mathd x} \ln \left( \Gamma \left( x \right) \right)$ is the
digamma function. The functions $\zeta_w \left( N ; s \right)$ have real zeros
at $s = - 1$ and $s = 0$, that is
\begin{equation}
  \lim_{s \rightarrow - 1} \zeta_w \left( N ; s \right) = \lim_{s \rightarrow
  0} \zeta_w \left( N ; s \right) = 0
\end{equation}
One possible idea is that the functions $\zeta_w \left( N ; s \right)$ can be
orthonormalized over the interval $s = - 1 \ldots 0$ via the Gram-Schmidt
process{\cite{matrix-computations}} and that the result might possibly shed
some light on the zeroes of $\zeta \left( s \right)$. Let the logarithmic
integral be defined
\begin{equation}
  \begin{array}{ll}
    \tmop{Li} \left( x \right) & =
  \end{array} \int_0^{\ln \left( x \right)} \frac{e^y - 1}{y} \mathd y + \ln
  \left( \ln \left( x \right) \right) + \gamma
\end{equation}
where $\gamma = 0.57721 \ldots$ is Euler's constant, then the normalization
factors are given by the integral
\begin{equation}
  \begin{array}{ll}
    \int_{- 1}^0 \zeta_w \left( N ; s \right) \mathd s & = \int_{- 1}^0
    \sum_{n = 1}^N \frac{n (n + 1)^{- s} - n^{1 - s} + s n^{- s}}{s - 1}
    \mathd s\\
    & = 1 + \frac{N}{N + 1} \left( \tmop{Li} \left( N + 1 \right) - \tmop{Li}
    \left( \left( N + 1 \right)^2 \right) \right) + \sum_{n = 1}^{N - 1}
    \frac{n}{\ln \left( n + 1 \right)}
  \end{array} \label{wmzi}
\end{equation}
\begin{figure}[h]
  \resizebox{5in}{3.25in}{\includegraphics{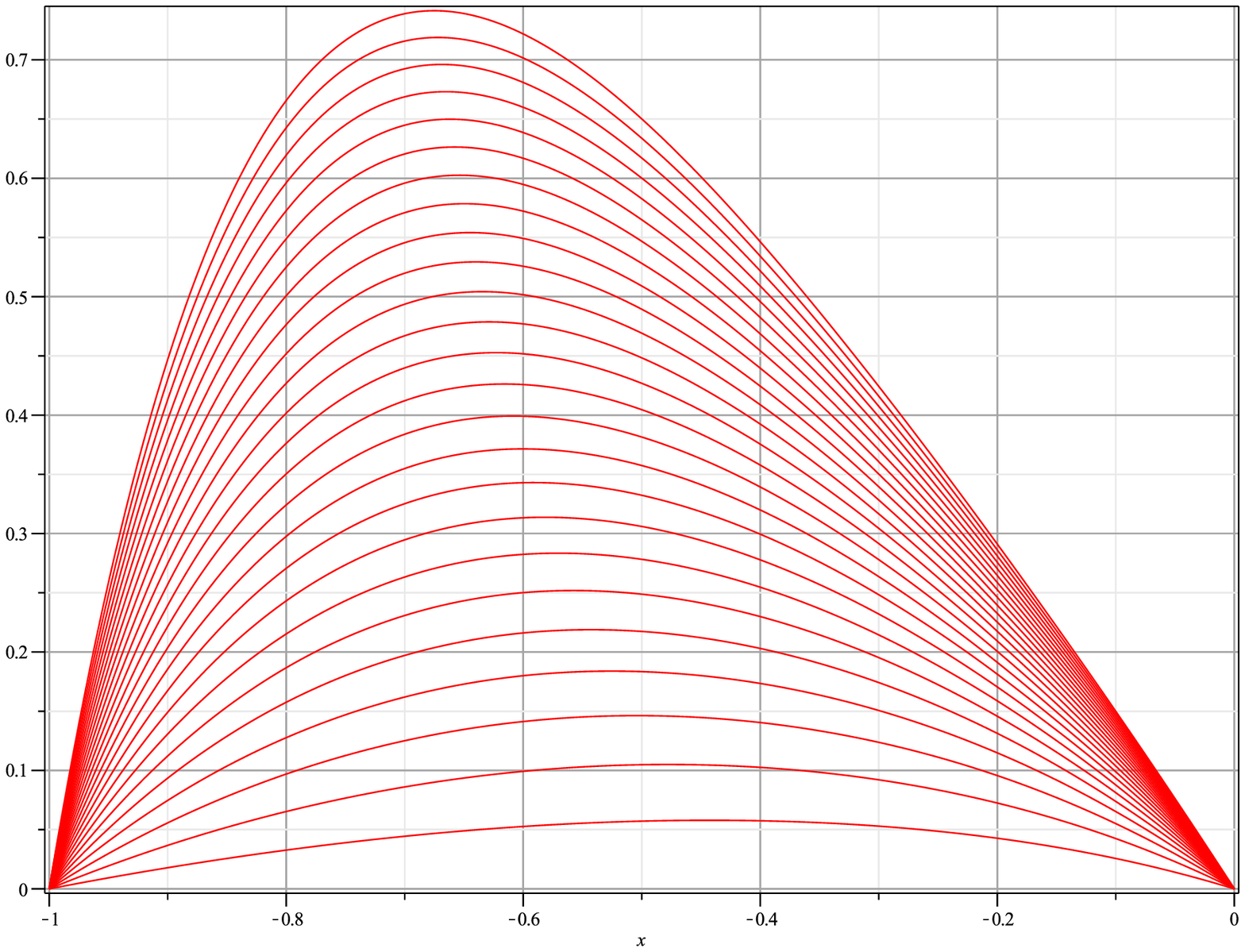}}
  \caption{$\left\{ \zeta_w \left( N ; s \right) : s = - 1 \ldots 0, N = 1
  \ldots 25 \right\}$}
\end{figure}

The following table lists the values of $\zeta_w \left( N ; 1 - n \right)$ for
$n = 2 \ldots 12$.
\[ \left[ \begin{array}{c}
     0\\
     - \frac{1}{6}  \hspace{0.25em} N - \frac{1}{6}  \hspace{0.25em} N^2\\
     - \frac{1}{4} N - \frac{1}{2} N^2 - \frac{1}{4}  \hspace{0.25em} N^3\\
     - \frac{7}{30} N - \frac{4}{5} N^2 - \frac{13}{15}  \hspace{0.25em} N^3 -
     \frac{3}{10} N^4\\
     - \frac{1}{6} \hspace{0.25em} N - \frac{11}{12}  \hspace{0.25em} N^2 -
     \frac{5}{3} N^3 - \frac{5}{4} N^4 - \frac{1}{3} N^5\\
     - \frac{5}{42}  \hspace{0.25em} N - \frac{6}{7} N^2 - \frac{97}{42} 
     \hspace{0.25em} N^3 - \frac{20}{7}  \hspace{0.25em} N^4 - \frac{23}{14} 
     \hspace{0.25em} N^5 - \frac{5}{14}  \hspace{0.25em} N^6\\
     - \frac{1}{8} N - \frac{19}{24}  \hspace{0.25em} N^2 - \frac{21}{8} 
     \hspace{0.25em} N^3 - \frac{14}{3} \hspace{0.25em} N^4 - \frac{35}{8} 
     \hspace{0.25em} N^5 - \frac{49}{24}  \hspace{0.25em} N^6 - \frac{3}{8}
     N^7\\
     - \frac{13}{90}  \hspace{0.25em} N - \frac{8}{9}  \hspace{0.25em} N^2 -
     \frac{26}{9}  \hspace{0.25em} N^3 - \frac{56}{9}  \hspace{0.25em} N^4 -
     \frac{371}{45}  \hspace{0.25em} N^5 - \frac{56}{9}  \hspace{0.25em} N^6 -
     \frac{22}{9}  \hspace{0.25em} N^7 - \frac{7}{18}  \hspace{0.25em} N^8\\
     - \frac{1}{10} N - \frac{21}{20}  \hspace{0.25em} N^2 - \frac{18}{5} 
     \hspace{0.25em} N^3 - \frac{79}{10}  \hspace{0.25em} N^4 - \frac{63}{5} 
     \hspace{0.25em} N^5 - \frac{133}{10}  \hspace{0.25em} N^6 - \frac{42}{5} 
     \hspace{0.25em} N^7 - \frac{57}{20}  \hspace{0.25em} N^8 - \frac{2}{5}
     N^9\\
     - \frac{1}{66}  \hspace{0.25em} N - \frac{10}{11}  \hspace{0.25em} N^2 -
     \frac{101}{22}  \hspace{0.25em} N^3 - \frac{120}{11}  \hspace{0.25em} N^4
     - \frac{199}{11}  \hspace{0.25em} N^5 - \frac{252}{11}  \hspace{0.25em}
     N^6 - \frac{221}{11}  \hspace{0.25em} N^7 - \frac{120}{11} 
     \hspace{0.25em} N^8 - \frac{215}{66}  \hspace{0.25em} N^9 - \frac{9}{22} 
     \hspace{0.25em} N^{10}\\
     - \frac{1}{12} N - \frac{1}{2} N^2 - \frac{55}{12}  \hspace{0.25em} N^3 -
     \frac{121}{8}  N^4 - \frac{55}{2} N^5 - \frac{110}{3}  \hspace{0.25em}
     N^6 - \frac{77}{2}  \hspace{0.25em} N^7 - \frac{231}{8}  \hspace{0.25em}
     N^8 - \frac{55}{4}  \hspace{0.25em} N^9 - \frac{11}{3} N^{10} -
     \frac{5}{12}  \hspace{0.25em} N^{11}
   \end{array} \right] \]

\subsubsection{Integrating Over the Critical Strip}

There is a formula similiar to (\ref{wmzi}) which gives the integral of
$\zeta_w \left( N ; s \right)$ over the critical strip $0 \leqslant \tmop{Re}
\left( s \right) \leqslant 1$.
\begin{equation}
  \begin{array}{ll}
    \int_0^1 \zeta_w \left( N ; c + i s \right) \mathd c & = 1 + \frac{N}{N +
    1} \left( \tmop{Ei}_1 \left( i s \ln \left( N + 1 \right) - \ln \left( N +
    1 \right) \right) - \tmop{Ei}_1 \left( i s \ln \left( N + 1 \right)
    \right) \right) + \sum_{n = 1}^{N - 1} \frac{n \left( n + 1 \right)^{- i
    s}}{\left( n + 1 \right) \ln \left( n + 1 \right)}
  \end{array}
\end{equation}
where $\tmop{Ei}_1 \left( t \right)$ is the exponential integral defined by
\begin{equation}
  \tmop{Ei}_1 \left( t \right) = t \int_0^1 \int_0^1 e^{- t x y} \mathd y
  \mathd x - \gamma - \ln \left( t \right)
\end{equation}
The contribution from the $\tmop{Ei}$ term vanishes as $s \rightarrow \infty$,
that is
\begin{equation}
  \lim_{s \rightarrow \infty} \frac{N}{N + 1} \left( \tmop{Ei}_1 \left( i s
  \ln \left( N + 1 \right) - \ln \left( N + 1 \right) \right) - \tmop{Ei}_1
  \left( i s \ln \left( N + 1 \right) \right) \right) = 0
\end{equation}

\subsubsection{The Reflection Formula}

There is a reflection equation for the finite-sum approximation $\zeta_{w_{}}
(N ; s)$ which is similiar to the well-known formula $\zeta \left( 1 - s
\right) = \chi \left( s \right) \zeta \left( s \right)$ with $\chi \left( s
\right) = 2 \left( 2 \pi \right)^{- s} \cos \left( \frac{\pi s}{2} \right)
\Gamma \left( s \right)$. The solution to
\begin{equation}
  \zeta_w \left( N ; 1 - s \right) = \chi \left( N ; s \right) \zeta_w \left(
  N ; s \right)
\end{equation}
is given by the expression
\begin{equation}
  \begin{array}{ll}
    \chi \left( N ; s \right) & = \frac{\zeta_w \left( N ; 1 - s
    \right)}{\zeta_w \left( N ; s \right)}\\
    & = \frac{\sum_{n = 1}^N - \frac{- n^s + \left( n + 1 \right)^{s - 1} n +
    n^{s - 1} - n^{s - 1} s}{s}}{\sum_{n = 1}^N \frac{- n^{1 - s} + \left( n +
    1 \right)^{- s} n + n^{- s} s}{s - 1}}\\
    & = - \frac{\left( s - 1 \right) \sum_{n = 1}^N - n^s + \left( n + 1
    \right)^{s - 1} n + n^{s - 1} - n^{s - 1} s}{s \sum_{n = 1}^N - n^{1 - s}
    + \left( n + 1 \right)^{- s} n + n^{- s} s}
  \end{array}
\end{equation}
which satisfies
\begin{equation}
  \chi \left( N ; 1 - s \right) = \chi \left( N ; s \right)^{- 1}
\end{equation}
The functions $\chi \left( N ; s \right)$, indexed by $N$, have singularities
at $s = 0$. Let
\begin{equation}
  \begin{array}{ll}
    a \left( N \right) & = \sum_{n = 1}^N n \left( \ln \left( n + 1 \right) -
    \ln \left( n \right) \right)\\
    b \left( N \right) & = \sum_{n = 1}^N \frac{\ln \left( n \right) n^2 - \ln
    \left( n + 1 \right) n^2 - \ln \left( n \right)}{n \left( n + 1 \right)}\\
    c \left( N \right) & = \frac{1}{2} \sum_{n = 1}^N n \left( \ln \left( n +
    1 \right)^2 - \ln \left( n \right)^2 \right)
  \end{array}
\end{equation}
then the residue at the singular point $s = 0$ is given by the expression
\begin{equation}
  \begin{array}{ll}
    \underset{s = 0}{\tmop{Res}} (\chi_{_{}} (N ; s)) & = - \underset{s =
    1}{\tmop{Res}} (\chi_{_{}} (N ; s)^{- 1})\\
    & = \frac{1 + \gamma + \Psi \left( n + 2 \right) - \frac{2}{N + 1} + b
    \left( N \right) - \frac{N \left( \ln \left( \Gamma \left( N + 1 \right)
    \right) - c \left( N \right) \right)}{\left( N - a \left( N \right)
    \right) \left( N + 1 \right)}}{a \left( N \right) - N}\\
    & = \frac{1 + \gamma + \Psi \left( n + 2 \right) - \frac{2}{N + 1} +
    \sum_{n = 1}^N \frac{\ln \left( n \right) n^2 - \ln \left( n + 1 \right)
    n^2 - \ln \left( n \right)}{n \left( n + 1 \right)} - \frac{N \left( \ln
    \left( \Gamma \left( N + 1 \right) \right) - \frac{1}{2} \sum_{n = 1}^N n
    \left( \ln \left( n + 1 \right)^2 - \ln \left( n \right)^2 \right)
    \right)}{\left( N - \sum_{n = 1}^N n \left( \ln \left( n + 1 \right) - \ln
    \left( n \right) \right) \right) \left( N + 1 \right)}}{\left( \sum_{n =
    1}^N n \left( \ln \left( n + 1 \right) - \ln \left( n \right) \right)
    \right) - N}
  \end{array}
\end{equation}
which has the limit
\begin{equation}
  \lim_{N \rightarrow \infty} \underset{s = 0}{\tmop{Res}} (\chi_{_{}} (N ;
  s)) = 1
\end{equation}
We also have the residue of the reciprocal at $s = 2$
\begin{equation}
  \begin{array}{ll}
    \underset{s = 2}{\tmop{Res}} (\chi_{_{}} (N ; s)^{- 1}) & = \frac{\frac{2
    N}{\left( N + 1 \right)^2} - 2 \Psi \left( 1, N + 1 \right) + 2 \zeta
    \left( 2 \right)}{\frac{\left( N + 1 \right)^2}{2} - \frac{N}{2} -
    \frac{1}{2} - \sum_{n = 1}^N n \left( \ln \left( n + 1 \right) + \ln
    \left( n + 1 \right) n - \ln \left( n \right) - n \ln \left( n \right)
    \right)}
  \end{array}
\end{equation}
which vanishes as $N$ tends to infinity
\begin{equation}
  \lim_{N \rightarrow \infty} \underset{s = 2}{\tmop{Res}} (\chi_{_{}} (N ;
  s)^{- 1}) = 0
\end{equation}
As can be seen in the figures below, the residue at $s = 0$ changes sign from
negative to positive between the values of $N = 176$ and $N = 177$.

\begin{figure}[h]
  \resizebox{5in}{3in}{\includegraphics{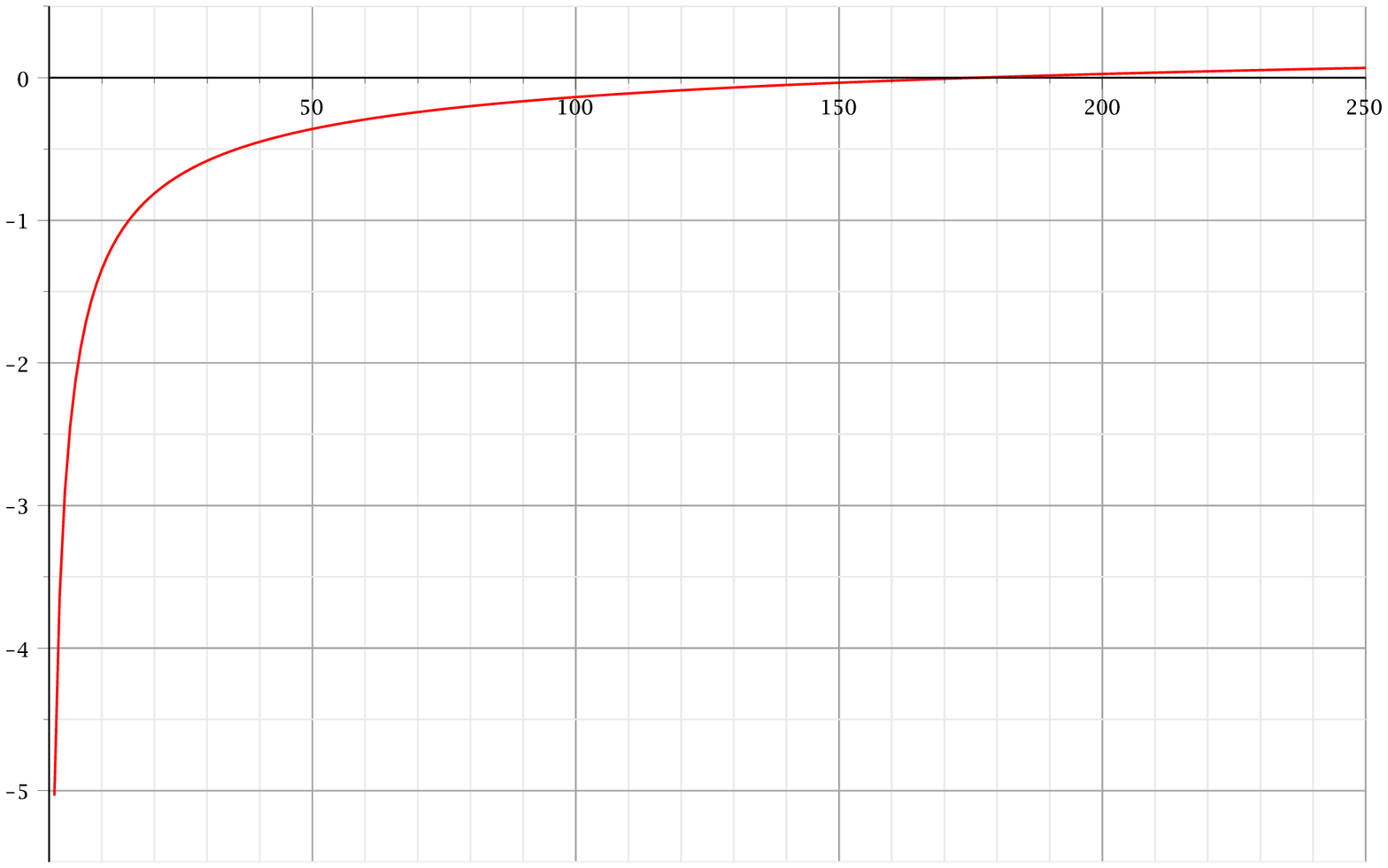}}
  \caption{$\left\{ \underset{s = 0}{\tmop{Res}} (\chi_{_{}} (N ; s)) : N = 1
  \ldots 250 \right\}$}
\end{figure}

\begin{figure}[h]
  \resizebox{5in}{3in}{\includegraphics{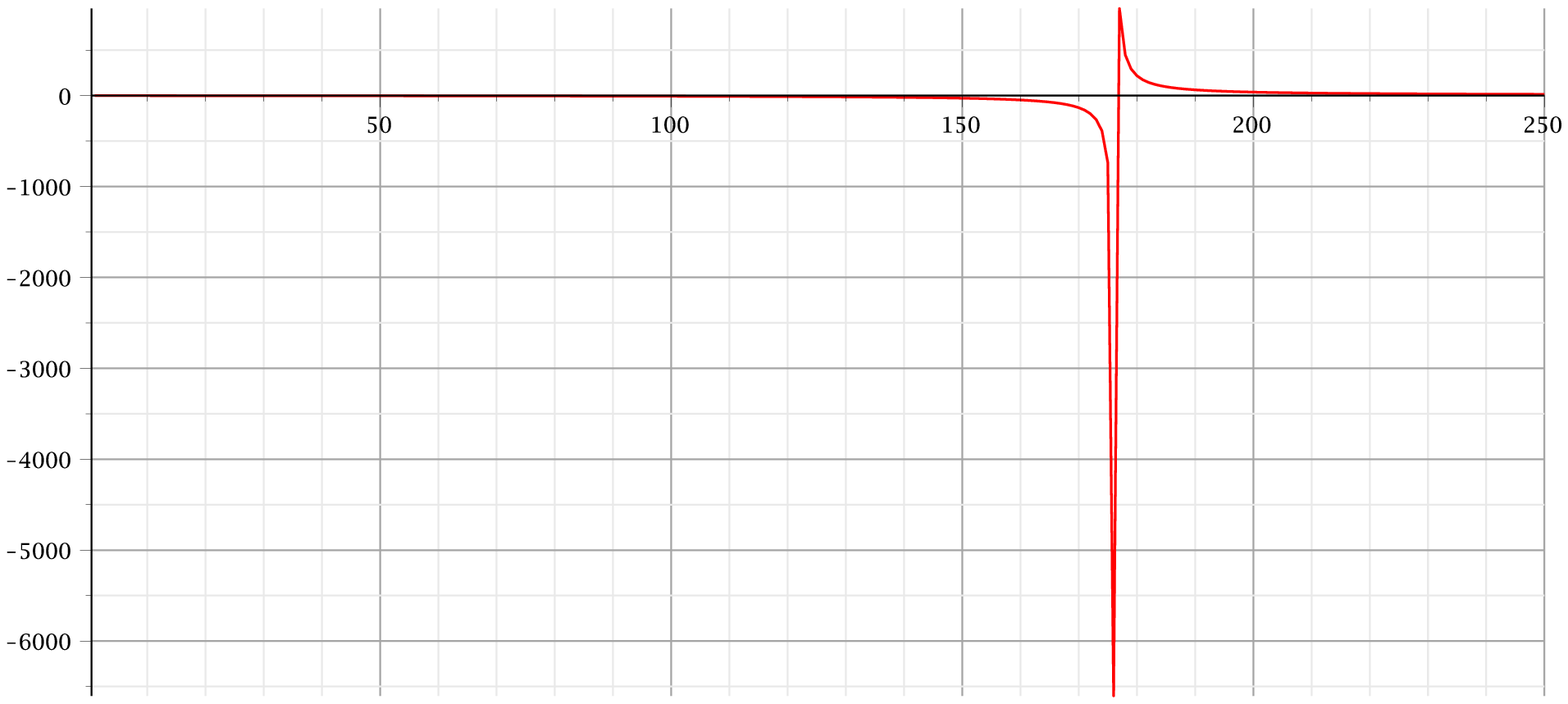}}
  \caption{$\left\{ \underset{s = 0}{\tmop{Res}} (\chi_{_{}} (N ; s))^{- 1} :
  N = 1 \ldots 250 \right\}$}
\end{figure}

For any positive integer N, we have the limits
\begin{equation}
  \begin{array}{ll}
    \lim_{s \rightarrow 0} \chi \left( N ; s \right) & = \infty\\
    \lim_{s \rightarrow 0} \frac{\mathd^n}{\mathd s^n} \chi \left( N ; s
    \right) & = \infty\\
    \lim_{s \rightarrow \frac{1}{2}} \chi \left( N ; s \right) & = 1\\
    \lim_{s \rightarrow 1} \chi \left( N ; s \right) & = 0\\
    \lim_{s \rightarrow 2} \chi \left( N ; s \right) & = 0\\
    \lim_{s \rightarrow 1} \frac{\mathd}{\mathd s} \chi \left( N ; s \right) &
    = 0
  \end{array}
\end{equation}
The line $\tmop{Re} \left( s \right) = \frac{1}{2}$ has a constant modulus
\begin{equation}
  \left| \chi \left( N ; \frac{1}{2} + i s \right) \right| = 1
\end{equation}
There is also the complex conjugate symmetry
\begin{equation}
  \chi \left( N ; x + i y \right) = \overline{\chi \left( N ; x - i y \right)}
\end{equation}
If $s = n \in \mathbbm{N}^{\ast}$ is a positive integer then $\chi \left( N ;
n \right)$ can be written as
\begin{equation}
  \begin{array}{cc}
    \chi \left( N ; n \right) & = \frac{\zeta_w \left( N ; 1 - n
    \right)}{\zeta_w \left( N ; n \right)}\\
    & = \frac{\sum_{m = 1}^N - \sum_{k = 1}^{n - 2} \frac{m^k}{n} \binom{n -
    1}{k - 1}}{\frac{N}{\left( n - 1 \right) \left( N + 1 \right)^n} -
    \frac{\cos \left( \pi n \right) \Psi \left( n - 1, N + 1 \right)}{\Gamma
    \left( n \right)} + \zeta \left( n \right)}\\
    & = \frac{- \sum_{m = 1}^N \frac{1}{n} \left( \left( n - 1 \right) m^{n -
    1} + m^n - \left( m + 1 \right)^{n - 1} m \right)}{\frac{N}{\left( n - 1
    \right) \left( N + 1 \right)^n} - \frac{\cos \left( \pi n \right) \Psi
    \left( n - 1, N + 1 \right)}{\Gamma \left( n \right)} + \zeta \left( n
    \right)}
  \end{array}
\end{equation}
where $\binom{n - 1}{k - 1}$ is of course a binomial. The Bernoulli
numbers{\cite{bernoilli-euler-maclaurin}} make an appearance since
\begin{equation}
  \begin{array}{ll}
    \chi \left( N ; 2 n \right) \zeta_{w_{}} (N ; 2 n) & = B_{2 n}  \left( N +
    1 \right)^2 \frac{\left( 2 n + 1 \right)}{2} + \ldots
  \end{array}
\end{equation}
The denominator of $\chi \left( N ; n \right)$ has the limits
\begin{equation}
  \begin{array}{cl}
    \lim_{N \rightarrow \infty} \zeta_w \left( N ; n \right) & = \zeta \left(
    n \right)\\
    \lim_{n \rightarrow \infty} \zeta_w \left( N ; n \right) & = 1
  \end{array}
\end{equation}
Another interesting formula gives the limit at $s = 1$ of the quotient of
successive functions
\begin{equation}
  \begin{array}{cl}
    \lim_{s = 1} \frac{\chi \left( N + 1 ; s \right)}{\chi \left( N ; s
    \right)} & = \frac{\left( N + 2 \right) N \left( N + 1 - a \left( N + 1
    \right) \right) }{\left( N + 1 \right)^2  \left( N - a \left( N \right)
    \right)}\\
    & = \frac{\left( N + 2 \right) N \left( N + 1 - \sum_{n = 1}^{N + 1} n
    \left( \ln \left( n + 1 \right) - \ln \left( n \right) \right) \right)
    }{\left( N + 1 \right)^2  \left( N - \sum_{n = 1}^N n \left( \ln \left( n
    + 1 \right) - \ln \left( n \right) \right) \right)}
  \end{array}
\end{equation}

\begin{figure}[h]
  \resizebox{5in}{3.25in}{\includegraphics{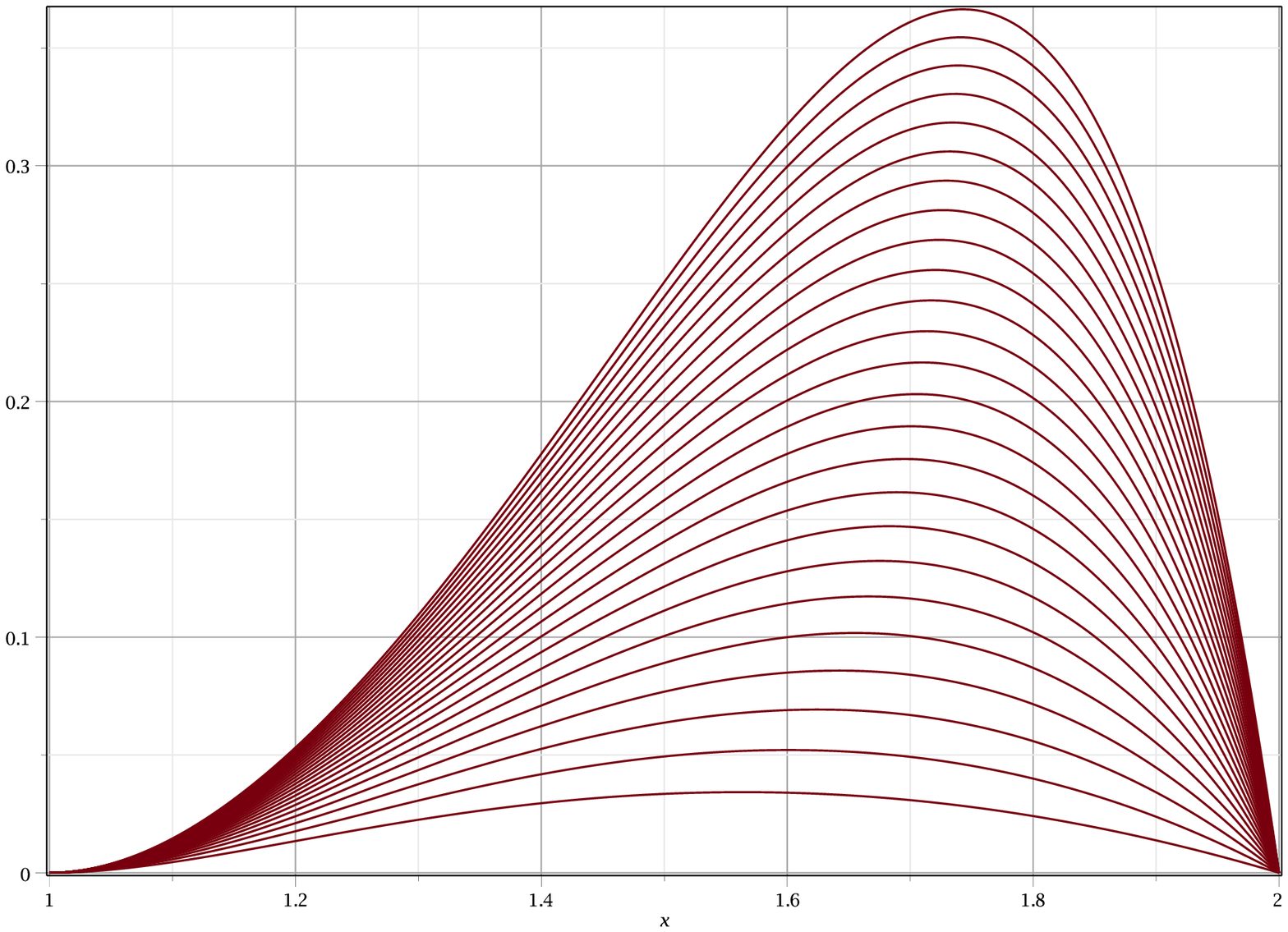}}

  \caption{$\left\{ \chi \left( N ; s \right) : s = 1 \ldots 2, N = 1 \ldots
  25 \right\}$}
\end{figure}

\begin{figure}[h]
  \resizebox{5in}{3.25in}{\includegraphics{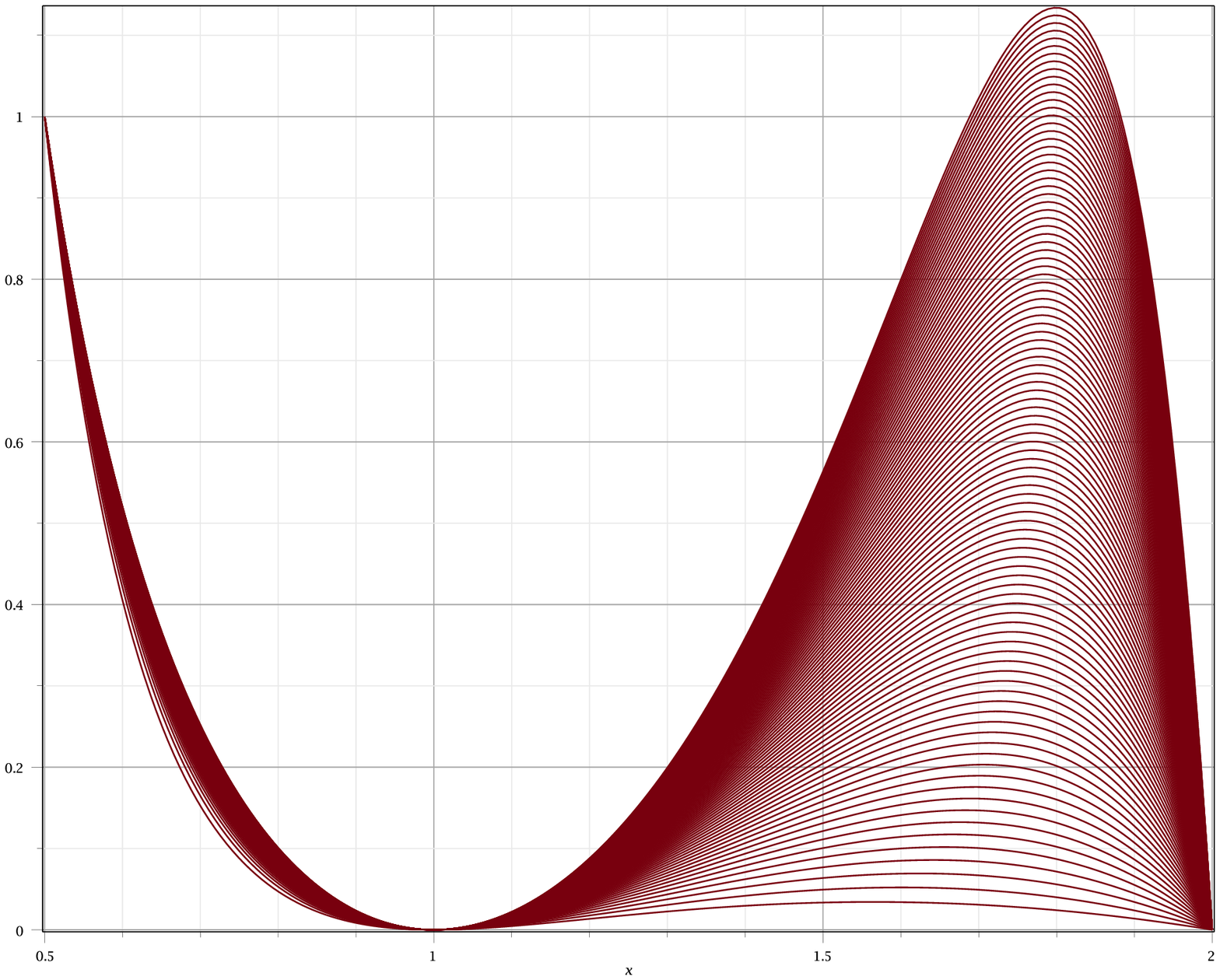}}

  \caption{$\left\{ \chi \left( N ; s \right) : s = \frac{1}{2} \ldots 2, N =
  1 \ldots 100 \right\}$}
\end{figure}

Let
\begin{equation}
  \nu \left( s \right) = \chi \left( \infty ; s \right) = \frac{\zeta \left( 1
  - s \right)}{\zeta \left( s \right)}
\end{equation}
Then the residue at the even negative integers \ is
\begin{equation}
  \underset{s = - n}{\tmop{Res}} (\nu_{_{}} (s)^{}) = \left\{
  \begin{array}{ll}
    \frac{\zeta \left( 1 - n \right)}{\frac{\mathd}{\mathd s} \zeta \left( s
    \right)_{} |_{s = - n}} & n \tmop{even}\\
    0 & n \tmop{odd}
  \end{array} \right.
\end{equation}


\begin{thebibliography}{1}
  \bibitem[1]{bernoilli-euler-maclaurin}G~Arfken.
  {\newblock}\tmtextit{Mathematical Methods for Physicists, 3rd ed.}, chapter
  5.9, Bernoulli Numbers, Euler-Maclaurin Formula., pages 327--338.
  {\newblock}Academic Press, 1985.
  
  \bibitem[2]{ithsm}Stephen Crowley. {\newblock}Integral Transforms of the Harmonic Sawtooth Map, The Riemann Zeta Function, Fractal Strings, and a Finite Reflection Formula {\newblock}http://arxiv.org/abs/1210.5652, October 2012.
  
  \bibitem[3]{two-new-zeta-constants}Stephen Crowley. {\newblock}Two new zeta
  constants: Fractal string, continued fraction, and hypergeometric aspects of
  the riemann zeta function. {\newblock}http://arxiv.org/abs/1207.1126, July
  2012.
  
  \bibitem[4]{matrix-computations}G.~H. Golub and C.~F. Van Loan.
  {\newblock}\tmtextit{Matrix Computations}. {\newblock}Johns Hopkins, 3
  edition, 1996.
\end{thebibliography}
\end{document}